\documentclass[12pt]{article}
\usepackage{graphicx}\usepackage{bm}
\usepackage{amsmath,amssymb,amsfonts,enumerate}
\def\Ex{\mathsf{E}}
\def\Var{\mathsf{V}}

\newtheorem{theorem}{Theorem}[section]

\newtheorem{prop}[theorem]{Proposition}
\newtheorem{cor}[theorem]{Corollary}

\begin{document}

\title{A dynamic network in a dynamic population: asymptotic properties}
\author{Tom Britton\thanks{Department of Mathematics, Stockholm
    University, SE-106 91 Stockholm, Sweden, tomb@math.su.se}, Mathias
  Lindholm\thanks{Uppsala University, Uppsala, Sweden. Current address: AFA Insurance, Stockholm, Sweden, Mathias.Lindholm@afaforsakring.se} and Tatyana
  Turova\thanks{Mathematical Center, Lund University University,
    Sweden, tatyana@maths.lth.se}}

\date{\today}
\maketitle

\begin{abstract}
\noindent We derive asymptotic properties for a
stochastic dynamic network model in a stochastic dynamic
population. In the model, nodes give birth to new nodes until they
die, each node being equipped with a social index given at
birth. During the life of a node it creates edges to other nodes,
nodes with high social index at higher rate, and edges disappear
randomly in time. For this model we derive criterion for when a giant
connected component exists after the process has evolved for a long
period of time, assuming the node population grows to infinity. We also obtain an explicit expression for the degree
correlation $\rho$ (of neighbouring nodes) which shows that $\rho$  is always
positive irrespective of parameter values in one of the two treated
submodels, and may be either positive or negative in the other model, depending on
the parameters.

\end{abstract}


\vskip.5cm
2010 Mathematics classification subjects. Primary: 92D30, Secondary: 60J80.

\section{Introduction}\label{sec-intro}

The models of dynamical
graphs are defined  by the rules of
attachment and deletion of vertices and edges chosen to fit particular
processes in nature. Intensive study in this area began in the end of
the 20th century when a number of models were formulated primarily in
the physics literature (see, e.g., Barab{\'a}si {\it et al.}\ (1999),
Callaway {\it et al.}\ (2001), Malyshev (1998)).
The important paper by Bollob{\'a}s,   Janson and
Riordan (2007) provided a unified
approach based on branching processes to many known  models  of random networks
(see the reference list in Bollob{\'a}s {\it et al.}\ (2007))

In the present paper we continue the study of a Markovian model
describing a random time-dynamic network in a random time-dynamic
population. The model was originally defined by Britton and Lindholm
(2010) extending an earlier model of Turova (2003), which in turn was
derived from a general model of Malyshev (1998). The population
process is a Markovian linear birth-and-death process in which individuals are
assigned random i.i.d.\ social indices. Given the population process, a
Markovian network is defined where edges between individuals appear
and disappear in such a way that ''social'' individuals tend to have
more neighbours. Notice also, that it was
shown by Turova (2002) that
the model without social index and
also without deletion of nodes includes as a subcase yet another model studied
by Callaway {\it et al.}\ (2001).

We study asymptotic properties of the network: what
properties will the network (and population) have after having evolved
for a long time assuming the size of the population has grown
large.
One can check that  
a snap-shot of the limiting network of
Britton and Lindholm
(2010)
falls  into the general class of
inhomogeneous random graphs introduced by Bollob{\'a}s, Janson and
Riordan (2007).
In particular, one can consult Britton and Lindholm
(2010) as well as Bollob{\'a}s, Janson and
Riordan (2007) 
for the age, the  type and  the degree distribution in this model.

In the
present paper we shall make
essential use of the theory of Bollob{\'a}s, Janson and
Riordan (2007) to study the phase
transition in the model by Britton and Lindholm (2010). 
More precisely, we will determine the
critical values which separate the set of parameters under which
the model with a high probability has a giant connected component (i.e.\ of order of the
entire network), and the area of parameters which do not
produce a giant component. This extends 
earlier results of
Turova (2007) for the model without social index.

Besides global properties of the network, e.g.\ phase transition, we
shall study here a characteristic, which tells more about the local
structure, namely the degree correlation, or the mixing coefficient as
it was introduced by Callaway {\it et al.}\ (2001). Inspiring numerical results on the
mixing coefficient in different empirical networks were presented and
analyzed already by Newman (2002). Recently Bollob{\'a}s,   Janson and
Riordan (2011) derived formulae for the mixing coefficient
in terms of small subgraphs counts for a rather general graph model,
which includes, in particular, the inhomogeneous graphs.

According to
empirical results provided by Newman (2002)
real life networks appear to be of two classes: assortative, when
the mixing  coefficient is positive, these are primarily the social
networks, or disassortative, when
the mixing  coefficient is negative, which is a feature of
technological or biological networks. Still, most most of the models
tend to have positive degree correlation/mixing coefficient. However,
Newman (2002) could construct  an example of a network model where
the mixing  coefficient changes sign
depending on the parameters of the model. Lately
Bollob{\'a}s,   Janson and
Riordan (2011) gave another recipe to construct a network with a negative
mixing coefficient. Therefore it is a challenge to look for an example of
a somewhat ``naturally grown'' random network which possesses the
property of enabling both positive and negative mixing coefficient
depending on the parameters. As it
is pointed out by Newman (2002) the property of assortativity affects
qualitatively the sharpness of the phase transition with respect to the
size of the giant component: the phase transition is sharper for the
dissasortative networks.

Here we provide an example  of a dynamic (growing) social network with a
degree correlation (mixing coefficient) that may be either negative or
positive depending on model parameters. We derive an explicit formula for the degree
correlation, showing the
dependence of the sign of the  mixing coefficient on the parameters.
Notice, that the formula for
the  mixing coefficient for
the model of Callaway {\it et al.}\ (2001) derived by Newman (2002),
is a subcase of the formula we prove here.

Although our model is a subcase of the general model studied in
Bollob{\'a}s,   Janson and
Riordan (2011), and hence the general formula from this paper should be applicable
for our model as well, we find a somewhat more direct way to get a formula for the degree
correlation in our case. We use the invariant measures
for  the random walk associated with the graph. It is worth noting that
our method is not restricted to the particular model we study here,
but would work as well for a class of inhomogeneous random graphs. However
this will be a subject of a separate study.

\section{The Markovian random network in a Markovian dynamic
population}

Below we define the Markovian random network in a Markovian dynamic
population, originally defined in Britton and Lindholm (2010). There are two
versions of the model, referred to as the Uniform (U) version and the
Preferential (P) version. In what follows a network denotes a finite set of nodes (the population) together with undirected edges connecting pairs of nodes. Nodes that are directly connected by an edge are called neighbours (in sociological applications nodes correspond to individuals, edges to some type of friendship, and neighbours are referred to as friends). The model is dynamic in the double sense that nodes are born and may die, and the same applies to edges.

\subsection{The model and its two versions}\label{sec-model}

We first define the node population dynamics. Let $Y(t)$ denote the
number of nodes alive at time $t$, and assume that $Y(0)=1$. While
alive, each node gives birth to new nodes at the constant rate
$\lambda$ and each node lives for an exponentially distributed time
having mean $1/\mu$ (denoted ${\rm Exp}(\mu)$), so each node dies at the rate $\mu$. We assume that
$\lambda >\mu$ implying that the expected number of births during a
life-span is larger than 1. This means that the node population is
modelled by a Markovian super critical branching process. Additional
to this, each node $i$ is at birth given a random ``social index''
$S_i$ having distribution $F_S$ on ${\mathbb R}^+$, independent and
identically distributed for different nodes. We will throughout assume
that $S$ has finite mean $\Ex[S]:=\mu_S<\infty$.

There are two versions of the model for births and deaths of edges,
both being Markovian given the node population and their social
indices. In both versions, nodes are isolated at birth, i.e.\ have no
neighbours. During the life of a node $i$, having social index $S_i=s$,
it creates new neighbour edges at rate $\alpha s$, also the same for
both models. The difference between the two versions lie in how
neighbour nodes are selected. In the uniform (U) version the neighbour
node is chosen uniformly at random among all living nodes. In the
preferential (P) version of the model the neighbour node is instead chosen at random among all living nodes with probabilities proportional to their social index. Finally, in both versions
each edge is removed, independently of everything else, at the rate $\beta$. If a node dies all edges connected to the node in question are removed. We emphasize that a node gets new neighbours in two ways: it ''creates'' new neighbour edges itself but may also be selected as neighbour of another nodes that has created a neighbour edge.

\subsection{Comments on the model}\label{sec-comments}

The model has four parameters: the birth and death rates of nodes,
$\lambda$ and $\mu$ respectively, and the death rate $\beta$ of edges
and $\alpha$ which is related to the birth rate of new edges. Beside
these four parameters there is the distribution $F_S$ for the social
indices $\{S_i\}$ (in what follows we let $f_s$ denote the
corresponding density function).

The P-version of the model is inspired by the preferential attachment
model \cite{BA99}, still being different in that here the probability
of receiving edges is determined at birth whereas it is determined by
random events during life in the preferential attachment model. In
both models there is also an age-factor in the sense that older nodes
tend to have more edges.

Some submodels are worth mentioning. One is where there is no
node-heterogeneity and nodes have the same social index $S\equiv s$
(for example set to 1 without loss of generality). Loosely speaking,
allowing node-heterogeneity makes it possibly to achieve degree
distributions having heavy tails: if $F_S$ is heavy tailed there will
be some nodes with very high social indices that hence have a large
number of neighbours (cf.\ the next section). The case where $\mu=0$ and $S\equiv 1$ has been
studied by Turova (2002) and (2007) who derive more results for this submodel.

\subsection{Known results of the model}\label{sec-old-results}

The model allows for multiple edges and self-loops which usually does not make sense in applications. However, it was shown by Britton and Lindholm (2010) that the \emph{proportion} of such edges is asymptotically negligible as long as $\Ex[S]<\infty$ for the U-version, and as long as $\Ex[S^2]<\infty$ for the P-version of the model. As a consequence the network will then have identical properties if loops and multiple edges are ignored or not allowed.

In the rest of the paper we will only consider the case where the node population tends to infinity, i.e.\ we condition on the event $B:=\{ \omega ;Y(t)\to \infty \}$,  which has positive probability since it was assumed that $\lambda>\mu$.
The first two results follow from standard branching process theory:

\emph{Asymptotic population size}: As $t\to\infty$
\[
Y(t)\sim e^{t(\lambda - \mu)}.
\]

\emph{Stable age distribution}: As $t\to\infty$ the age $A$ of a
randomly selected living node will satisfy $A\sim {\rm Exp}(\lambda)$
having mean $1/\lambda$, see e.g.\ Section 3.4 in Haccou et al. (2005).

\emph{Social index distribution}: Since the social indices are defined to be i.i.d.\ having distribution $F_S$, and they only affect the occurrence of edges, the social index of a randomly selected living node will have the same social index distribution $F_S$.

\emph{Asymptotic degree distribution}: As $t\to\infty$ the degree
(i.e., the number of neighbours) $D$ of a randomly selected living individual has a mixed Poisson distribution in both versions of the model (Britton and Lindholm, 2010). For the U-version it is
\begin{equation}
D^{(U)}\sim {\rm MixPo}\left(
\frac{\alpha (S+\mu_S)\left(1-e^{-(\beta +\mu) A}\right)}{\beta +\mu}\right),\label{uncond-distr}
\end{equation}
where $A\sim {\rm Exp}(\lambda)$ and $S\sim F_S$ are independent. This simply means that,
 conditional on $A=a$ and $S=s$, the degree is Poisson distributed
 with mean parameter equal to $\alpha (s+\mu_S)\left(1-e^{-(\beta +\mu)
     a}\right)/(\beta +\mu)$. The mean and variance of this distribution
are given by
\begin{align*}
\Ex[D^{(U)}]&=\frac{2\alpha}{\lambda + \beta +\mu}\mu_S,\\
\Var[D^{(U)}] &=\frac{2\alpha}{\lambda+\beta +\mu}\mu_S \\
& +
\frac{4\lambda\alpha^2}{(\lambda + \beta +\mu)^2(\lambda+2(\beta +\mu))}\mu_S^2
+ \frac{2\alpha^2}{(\lambda+\beta +\mu)(\lambda+2(\beta +\mu))}\Var[S].
\end{align*}

For the P-version of the model the degree distribution satisfies
\begin{equation}
D^{(P)}\sim {\rm MixPo}\left(\frac{2\alpha S\left(1-e^{-(\beta +\mu) A}\right)}{\beta +\mu}\right), \label{mod-uncond-distr}
\end{equation}
where $A\sim {\rm Exp}(\lambda)$ and $S\sim F_S$ are independent. The mean and variance are given by
\begin{align*}
\Ex[D^{(U)}]&=\frac{2\alpha}{\lambda + \beta +\mu}\mu_S,\\
\Var[D^{(U)}] &=\frac{2\alpha}{\lambda+\beta +\mu}\mu_S \\
& +
\frac{4\lambda\alpha^2}{(\lambda + \beta +\mu)^2(\lambda+2(\beta +\mu))}\mu_S^2
+ \frac{8\alpha^2}{(\lambda+\beta +\mu)(\lambda+2(\beta +\mu))}\Var[S],
\end{align*}
i.e.\ the same mean degree but a larger variance. For both versions it
is seen that the variance of the degree distribution increases with the
variance of the social index distribution, so having a heavy-tailed
social index distribution $F_S$ implies that also the degree
distribution will be heavy-tailed.

\section{Results}\label{sec-results}

In what follows we derive some further properties of the model, still
assuming $t\to\infty$ and conditioning on that the node population
grows beyond all limits. We do this for the U-version and only present
the results for the P-version which is obtained in a similar way.

\subsection{Type-distribution of neighbour nodes}

Assume $t$ to be large and consider a randomly selected individual alive at this time
having age $a$ and social index $s$. We derive the distribution of
the number of $(a',s')$-neighbours  this $(a,s)$-node has. The number of neighbours of any type $(a',s')$ will be Poisson distributed, so we first compute the mean of the distribution and then divide by the total mean thus obtaining the type-distribution.

How many neighbours of type $(a',s')$ does the $(a,s)$-individual
have? Any such edge must have been created at some time $\tau\in
(t-a\wedge a', t)$ ('$\wedge$' denotes minimum: both individuals must have been born). At such
$\tau$ our $(a,s)$-individual creates such an edge at rate $\alpha s
\lambda e^{-\lambda (\tau-(t-a'))}f_S(s')$, because it creates edges
at rate $\alpha s$ and the second part denotes the fraction of nodes
having age $\tau -(t-a')$ (at $\tau$) and type $s'$. Similarly, it
receives edges at the rate $\alpha s' \lambda e^{-\lambda
  (\tau-(t-a'))}f_S(s')$ because there are $Y(\tau)\lambda e^{-\lambda
  (\tau-(t-a'))}f_S(s')$ having the ''correct'' age and type, each of
them creates edges at rate $\alpha s'$ and with probability
$1/Y(\tau)$ it reaches our $(a,s)$-individual. Any such edge created
at $\tau$ remains at $t$ if and only if both the $(a',s')$-individual
and the edge
survive until $t$ (we ave already
conditioned on that our $(a,s)$-individual lives at $t$). The
probability for this equals $e^{-(\beta+\mu)(t-\tau)}$. The total
expected number of edges with $(a',s')$-individuals hence equals
\begin{align}
m^{(U)}(a',s'|a,s)da'ds'&=da'ds'\int_{t-a\wedge a'}^t \alpha (s+s')\lambda e^{-\lambda
  (\tau-(t-a'))}f_S(s')e^{-(\beta+\mu)(t-\tau)} d\tau.
\\
&=\alpha (s+s')f_S(s')\lambda e^{-\lambda a'} \frac{e^{(\lambda -\beta -
    \mu)a\wedge a'}-1}{\lambda-\beta-\mu}da'ds'.\label{mean_spec_type}
\end{align}
 It is easy to check from this that the expected total number of neighbours of any type of the $(a,s)$-node, the integral of this quantity over all $a'$ and $s'$, equals
\begin{equation}
m^{(U)}(\cdot , \cdot |a,s)= \frac{\alpha (s+\mu_S)}{\beta +\mu}\left(1-e^{-(\beta +\mu)a}\right).\label{mean_total}
\end{equation}
This also confirms that the unconditional degree distribution has the mixed
Poisson distribution specified in Equation (\ref{uncond-distr}).

As a consequence, for the U-version of the model, the type
distribution of the neighbours of an individual of type $(a,s)$, which
equals the expected number of $(a',s')$-types divided by the total
expected number of types, becomes
\[
f^{(U)}(a',s'|a,s)=\frac{m^{(U)}(a',s'|a,s)}{m^{(U)}(\cdot , \cdot |a,s)}=\frac{(s+s')f_S(s')}{s+\mu_S} \frac{(\beta
   +\mu)\lambda e^{-\lambda a'} \left(e^{(\lambda -\beta-\mu)a\wedge
      a'}-1\right) }{(\lambda -\beta -\mu)(1-e^{-(\beta + \mu)a})}.
\]
We hence have conditional independence:
$f^{(U)}(a',s'|a,s)=f^{(U)}(a'|a)f^{(U)}(s'|s)$, where
\begin{align}
f^{(U)}(a'|a)&=\frac{(\beta
  +\mu)\lambda e^{-\lambda a'} \left(e^{(\lambda -\beta-\mu)a\wedge
      a'}-1\right) }{(\lambda -\beta -\mu)(1-e^{-(\beta + \mu)a})}\label{f_U(aa)}
,\\
f^{(U)}(s'|s) &= \frac{(s+s')f_S(s')}{s+\mu_S}.\label{f_U(ss)}
\end{align}
Admittedly, the notation here, and in what follows, is a bit sloppy in that $f^{(U)}(a'|a)$ and
$f^{(U)}(s'|s)$ really denote different density functions, as
indicated by the difference in arguments.

Using similar arguments, but for the P-version of the model, we get
$f^{(P)}(a',s'|a,s)=f^{(P)}(a'|a)f^{(P)}(s'|s)$ where
\begin{align}
f^{(P)}(a'|a)&=\frac{(\beta
  +\mu)\lambda e^{-\lambda a'} \left(e^{(\lambda -\beta-\mu)a\wedge
      a'}-1\right) }{(\lambda -\beta -\mu)(1-e^{-(\beta + \mu)a})} ,\label{f_P(aa)}\\
f^{(P)}(s'|s) &= \frac{s'f_S(s')}{\mu_S},\label{f_P(ss)}
\end{align}
i.e.\ the same density with respect to age but different with respect
to social index.

\subsection{Phase transition}\label{PT}

Here we shall apply the results of Bollob{\'a}s, Janson and Riordan (2007) to the introduced
above models in the limit  as $t\rightarrow
\infty$.
Therefore first we shall put our model into a general setup of the theory (and notation)
of inhomogeneous random graphs of Bollob{\'a}s, Janson and Riordan (2007).

Let us enumerate the individuals at time $t$ by $i=1, \ldots,
n:=Y(t)$. Let $X_i=(A_i, S_i)$ denote, correspondingly, the age and social
index  of individual $i$.
Notice, that this means that
this individual was born at time $t-A_i$ (and remains alive at time
$t$).
The random variables $X_i$ are independent for different $i$.
For a randomly selected $i$, as we argued above, the distribution of
the social index follows the original social index distribution and
$A_i$ is Exp$(\gamma)$; denote the corresponding measures by $\mu_1$
and $\mu_2$:
\begin{equation}\label{mu}
  \begin{array}{ll}
\mu_1(da)&=\lambda e^{-\lambda a} da,   \\
\mu_2(ds)&=f_S(s)ds.
\end{array}
\end{equation}
We shall call $X_i$ the type of individual $i$. Given the types
$X_i=(a_i, s_i)$ and $X_j=(a_j, s_j)$ of
individuals $i$ and $j$, we derive the probability $p_{i j}(n)$ that they are
connected. From above we know that the expected number of
$(a',s')$-individuals our $i$ is connected to equals
$m(a',s'|a_i,s_i)da'ds'$ defined by (\ref{mean_spec_type}) for the U-version
(and similarly for the P-version). In total, if $Y(t)=n$, the expected
total number of individuals of type $(a', s')$ equals
$n\lambda e^{-\lambda a'}da'f_S(s')ds'$. This reasoning implies that
the probability that $i$ is connected to $j$ (of type $(a_j,s_j)$)
equals the ratio of these two expressions: 
\begin{equation}\label{p1}
p_{i j}(n)=
\frac{1}{n} \kappa_1(a_i,a_j) \ \kappa_2(s_i,s_j),
\end{equation}
where
\[
\kappa_1(a_i,a_j)=\frac{e^{(\lambda -\beta -\mu) (a_i\wedge a_j)}-1}{\lambda -\beta
  -\mu}.
\]
and
\[\kappa_2(s_i,s_j)=\left\{
\begin{array}{ll}
  \kappa_2^U(s_i,s_j) =\alpha(s+s' )  &  \mbox{ in U-version, } \\ \\
\kappa_2^P(s_i,s_j)=2 \alpha ss'  &  \mbox{ in P-version. }
\end{array}
\right.
\]

Let
$\mu=\mu_1\times \mu_2$, which by (\ref{mu}) is a probability
measure on ${\cal S}={\bf R}_+^2$,
and let
\[
\kappa =
\kappa(x_i,x_j)=\kappa_1(a_i,a_j)
\kappa_2(s_i,s_j).
\]
Given the sequence $x_1, \ldots, x_{n}$, we let  $G^{\mathcal V}(n,\kappa)$
 be the random graph on $\{1, \ldots, n\}$, such that any
two vertices $i$ and $j$ are connected by an edge, independently of the others
vertices, with  probability given by (\ref{p1}), which we write as
\[
p_{i j}(n)= \min\{\kappa (x_i,x_j)/n,1\}.
\]
Hence, conditionally on $Y(t)=n$ the graph $G^{\mathcal V}(n,\kappa)$
describes our model. On the other hand, $G^{\mathcal V}(n,\kappa)$
with $\kappa$ as above satisfies the definition of inhomogeneous
random graph from Bollob{\'a}s, Janson and Riordan (2007), and thus we can apply some known results here.

Recall first the fundamental result on Phase Transitions in the
inhomogeneous random graphs.
Define
\[T_{ \kappa}f(x)=\int_S
  \kappa(x,y)f(y)
d\mu(y),
\]
and
\[\| T_{ \kappa}\|= \sup\{\| T_{\kappa}f\|_2:f\geq 0,
\|f\|_2\leq1\}.\]
 Then, by  Theorem
3.1 from \cite{BJR}, the largest connected component of $G^{\mathcal
  V}(n, \kappa)$, denoted $C_1 ( G^{\mathcal V}(n, \kappa) )$,
satisfies 
\[\frac{1}{n}
C_1( G^{\mathcal V}(n, \kappa) )
\stackrel{P}{\rightarrow} \rho_{\kappa}=\int_S
  \rho(x)
d\mu(x), \hbox{ as } n\to\infty,
\]
where  $\rho_{\kappa}>0$  if and only if $\| T_{\kappa}\|>1$.

This gives as well for our  model that
conditionally on the number  $Y(t)$ of individuals in the
population at time $t$, we have with a high probability approximately $\rho_{\kappa}Y(t)$ individuals connected; and
thus this holds
unconditionally as well.
The function  $\rho_{\kappa}(x)$ is described more precisely in the following
theorem.
\begin{theorem} [\cite{BJR}, Theorem 6.1]
  Suppose that  $\kappa$ is the kernel on $({ \cal S}, \mu)$,
  that
  $ \kappa \in
  L^1$, and
 \[
\int_{ \cal S}
\kappa(x,y)d\mu(y)<\infty
 \]
  for every $x\in S$.
Then $\rho_{\kappa} (x)$
 is the maximal
  solution to
\begin{equation}\label{+}
f(x) = 1- \exp\{-T_{ \kappa}f(x)\}.
\end{equation}
Furthermore:

$(i)$ If $\| T_{\kappa}\|\leq 1$ then $\rho_{\kappa}
(x)=0$ for every $x$, and (\ref{+})
has only the zero solution.

$(ii)$ If $1< \|T_{\kappa}\|\leq \infty$ then $\rho_{\kappa}
(x)>0$ on a set of a positive measure. If, in addition, $\kappa$
is irreducible, i.e., if
\[A\subseteq S \mbox{ and } \kappa=0 \ a.e. \ \mbox{ on } A\times(S\setminus A) \mbox{ implies }
\mu(A) \mbox{ or } \mu(S\setminus A)=0,\]
then $\rho_{\kappa}
(x)>0$
for a.e. $x$, and $\rho_{\kappa}
(x)$ is the only non-zero solution of (\ref{+}).
\end{theorem}

First we shall derive conditions for $\rho_{\kappa}
(x)>0$. By the cited results this amounts to computing $\|T_{\kappa}\|$.
Unfortunately, there is no general formula to simply compute the norm
$\|T_{\kappa}\|$.
Therefore
we shall use a criteria for the condition when
$\|T_{\kappa}\|>1$, which is what we need here.

\begin{prop}\label{P1} [\cite{BJR}, Proposition 17.2]
For $k\geq 1 $ let
\[\alpha (k):=\frac{1}{2} \int_{{ \cal S}^{k+1}}\kappa(x_0, x_1)\kappa(x_1, x_2)\ldots \kappa(x_{k-1},x_k) d \mu(x_0)\ldots d \mu(x_k).\]
Then $\|T_{\kappa}\|>1$ if and  only if
$\alpha(k) \rightarrow \infty$ as $k\to\infty$.
\end{prop}

It follows from definition (\ref{mu}) that for our model
\begin{equation}\label{al}
\alpha (k)= \frac{1}{2} \alpha_1 (k)\alpha_2 (k),
\end{equation}
where for $i=1,\ 2$,
\begin{equation}\label{ali}
\alpha_i(k)= \int_{{ \bf R}_+^{k+1}}\kappa_i(u_0, u_1)\kappa_i(u_1, u_2)\ldots \kappa_i(u_{k-1},u_k)d \mu_i(u_0)\ldots d\mu_i(u_k).
\end{equation}

First we shall compute $\alpha_2(k)$. We start with the U-version in the next
Proposition.
\begin{prop}\label{PY}
Introduce the matrix
\[A =\left(
\begin{array}{ll}
\Ex[S] & 1\\
\Ex[S^2] & \Ex[S]
\end{array}
\right).\]
Then for all $k\geq 3$
\begin{equation}\label{alk}
\alpha_2^U(k)=\alpha^k\left\{A^{k-1}\left(
\begin{array}{c}
2\Ex[S] \\
(\Ex[S] )^2+\Ex[S^2]
\end{array}
\right) \right\}_{11},
\end{equation}
where for a matrix $B $ we denote $\left\{B\right\}_{ij}$ the corresponding entry.
\end{prop}

\noindent
{\bf Proof.} Let $S_0, S_1, \ldots ,$ be $i.i.d.$\ copies of the random variable $S$. Then
by the definition (\ref{ali})
\begin{equation}\label{N1}
\alpha_2^U(k)=\alpha^k \Ex\left[\prod_{i=1}^k(S_{i-1}+S_i)\right]=:\alpha^k f(k),
\end{equation}
where
\[f(k)=\Ex\left[\prod_{i=1}^k(S_{i-1}+S_i)\right] .\]
Denote also
\[g(k)= \Ex\left[\left(\prod_{i=1}^k(S_{i-1}+S_i)\right)S_k\right] .\]
Then we recursively derive
\[f(k)= g(k-1)+f(k-1)\Ex[S],\]
and
\[g(k)= g(k-1)\Ex[S]+f(k-1)\Ex[S^2],\]
for $k>1$, with
\begin{equation}\label{al1}
f(1)=2 \Ex[S], \ \ \ \ \ g(1)=(\Ex[S])^2+\Ex[S^2].
\end{equation}
Hence, for all $k>1$
\[
\left(
\begin{array}{c}
f(k)\\
g(k)
\end{array}
\right) =A
\left(
\begin{array}{c}
f(k-1)\\
g(k-1)
\end{array}
\right)
=A^{k-1}\left(
\begin{array}{c}
f(1)\\
g(1)
\end{array}
\right),
 \]
and formula (\ref{alk}) follows from here,  (\ref{al1}) and
(\ref{N1}).$\hfill \Box$

\begin{cor}\label{C1}
  One has
\begin{equation}\label{N4}
  \lim_{k\rightarrow \infty}\left( \alpha_2^U(k)  \right)^{1/k}=
  \alpha \left(\Ex[S] +\sqrt{\Ex[S^2]}  \right).
  \end{equation}
\end{cor}

In the case of the P-version, when $\alpha_2=\alpha_2^P$ it is straightforward to
  derive from definition (\ref{ali}) that
  \begin{equation}\label{N2}
 \alpha_2^P(k)=(2\alpha)^k \left( \Ex[S]  \right)^{2}    \left(
   \Ex[S^2]  \right)^{k-2},
\end{equation}
which yields
\begin{equation}\label{N3}
   \lim_{k\rightarrow \infty}\left( \alpha_2^P(k)  \right)^{1/k}=
  2 \alpha  \Ex[S^2] .
\end{equation}

Define now
\begin{equation}\label{cr1}
c ^{cr} (\lambda , \mu + \beta ) := \sup \{ x>0 : \sum _{k=2}^{\infty}
x^k \  \alpha_1(k)
 < \infty \}.
\end{equation}
Then due to  the definitions  (\ref{al}) and (\ref{ali}) together with
asymptotics (\ref{N4}) and  (\ref{N3})  we have the following
criteria which holds for both cases, P-version and U-version.

\begin{cor}\label{C2}
If
\[\lim_{k\rightarrow \infty}\left( \alpha_2 (k)  \right)^{1/k} > c
^{cr} (\lambda , \mu + \beta )\]
then $\alpha (k) \rightarrow \infty $;
if
\[\lim_{k\rightarrow \infty}\left( \alpha_2 (k)  \right)^{1/k} < c
^{cr} (\lambda , \mu + \beta )\]
then $\alpha (k) < \infty $.
\end{cor}

The value $c
^{cr} (\lambda , \mu + \beta )$  is known
from \cite{T03}. Let us record this result here. Write first
\[\alpha_1(k) = \int_{{\bf R}_+^{k+1}} \kappa_1(a_0, a_1)\cdot \ldots
\cdot
\kappa_1(a_{k-1}, a_k)
\lambda e^{-\lambda a_0} da_0  \ldots\lambda e^{-\lambda a_k} da_k
\]
\[= \Ex\left[ \kappa_1\left(\frac{X_0}{\lambda}, \frac{X_1}{\lambda}\right)\cdot \ldots\cdot
\kappa_1\left(\frac{X_{k-1}}{\lambda}, \frac{X_k}{\lambda}\right)\right],
\]
where $X_0, \ldots , X_k$ are $i.i.d.$ random variables with a common
$\mbox{ Exp }  (1)$-distribution. Then
by the results of \cite{T07} (see \cite{T07} formula (1.7), or consult
\cite{BJR}) the value
$c ^{cr} (\lambda , \mu + \beta )$
is the smallest positive root of
\begin{equation}\label{FA}
H(x)=1+\sum _{n=1}^{\infty} (-1)^{n}
\left( \frac{x}{\mu+ \beta} \right)^n \,  \frac{1}{n!} \,
\left[  \, \prod_{l=1}^n
\frac{1}{1+(l-1)
\, \frac{\mu+\beta}{\lambda}}
\, \right].
\end{equation}

Finally, combining Corollary \ref{C2} and Proposition \ref{P1}
together with asymptotics (\ref{N4}) and  (\ref{N3})  we can derive,
apart from the critical case, the necessary and sufficient conditions
for the existence of the giant component.

\begin{cor}\label{C3}
  Denote
  \[{\kappa}^U={\kappa}_1{\kappa}_2^U, \ \
  {\kappa}^P={\kappa}_1{\kappa}_2^P,\]
and let 
\begin{equation}
R^U(\alpha, \beta, \lambda,\gamma, S) :=\frac{\alpha\left(\Ex[S]
    +\sqrt{\Ex[S^2]}  \right)}{c^{cr}(\lambda , \mu + \beta ) }
\quad \text{and}\quad
R^P(\alpha, \beta, \lambda,\gamma, S) :=\frac{2\alpha \Ex[S^2]}{c^{cr}(\lambda , \mu + \beta ) }.\label{R^U}
\end{equation}
Then there exists a giant component in the
  $U$-model if $R^U(\alpha, \beta, \lambda,\gamma, S)>1$,
  and there is no giant component if $R^U(\alpha, \beta, \lambda,\gamma, S)<1$. 
 
Similarly,   there exists a giant component in the
  $P$-model if $R^P(\alpha, \beta, \lambda,\gamma, S)>1$,
and there is no giant component if $R^P(\alpha, \beta, \lambda,\gamma, S)<1$.
\end{cor}

Note also that although we define the function $c^{cr} (\lambda , u )$ for
$u>0$, the model with $ \mu + \beta=0$ is treated in a very similar
way, in fact when also $S \equiv 1$, this model becomes the one
studied by Callaway {\it et al.}\ (2001). Then one defines $c^{cr} (\lambda , 0)$
as the critical parameter, above which there is a giant component. It
is known (see Callaway {\it et al.}\ (2001)) that $c^{cr}
(\lambda , 0)=1/4$, and moreover function $c^{cr} (\lambda , u )$ is
continuous at $u=0$ (see Bollob{\'a}s {\it et al.}\ (2007) and Turova (2007a)).

In general, there  is no closed form for the function $c^{cr}$, but
one can mention some qualitative properties of this function. First of
all, simply from the model it follows that  $c^{cr}(\lambda , \mu + \beta )$ is
decreasing in $\lambda$ and increasing in $\mu + \beta$, i.e., both in
$\mu$ and $\beta$.
For any $u>0$ denote $c(u)$ the smallest positive root of the rescaled
function (\ref{FA})
\[
H_u(x)=1+\sum _{n=1}^{\infty} (-1)^{n}
x^n \,  \frac{1}{n!} \,
\left[  \, \prod_{l=1}^n
\frac{1}{1+(l-1)
\, u}
\, \right].
\]
Then
\[c^{cr} (\lambda , \mu + \beta) = (\mu + \beta)c\left(\frac{\mu + \beta}{\lambda} \right).\]

\subsection{Degree correlation}\label{sec-deg-corr}

We now derive an expression for the degree correlation $\rho$ (also
known as ''mixing coefficient'') of the network
when it has grown and reached its stationary phase. That is, after a
long time $t$ we take a snap-shot of the network and compute the
degree correlation $\rho$, which is defined as the correlation of the
degrees of the two adjacent nodes of a randomly selected edge in the
network. We derive $\rho$ in two different ways. In the first method
we first pick a random \emph{node} in the snap-shot network and then
perform a random walk on this fixed network. When this random walk has
reached stationarity we consider the node in one step as the ''first''
node and then pick the ''second'' node randomly among the neighbouring
nodes, and let the edge between these two nodes be our randomly
selected edge. In the second method we pick our ''first'' node
randomly in a size-biased way in the snap-shot network (such that the
node is adjacent to a randomly selected \emph{edge}) and the
''second'' node randomly among the neighbours of this individual (the
size-biased distribution of a non-negative random variable $X$ with
density/pdf $f_X(x)$ and finite mean $\mu_X$ has density/pdf
$xf_X(x)/\mu_X$ -- below we denote such a random variable by $\tilde
X$). As will be seen, the two methods of picking neighbouring nodes
give the same result. The degree correlation is then obtained by
computing the degree correlation of the two neighbouring nodes.

\subsubsection{The type-distribution obtained using random walk method}

We start with the first method where we first pick a node at random among all living nodes. From the results of Section \ref{sec-old-results} this means the node has exponential age distribution and
independent social index distribution given by the original
distribution. More precisely the type distribution of a randomly
selected node equals $\lambda e^{-\lambda a}da f_S(s)ds$. The neighbours of this node, conditional on its type, all have
distribution $f(a',s'|a,s)$, defined in Section 3.1, where its form
depends on whether we consider the U-version or the P-version of the
model. Iterating this procedure gives the type-distribution of nodes
after additional steps in the random walk. Eventually this random walk
reaches stationarity and the distribution $f_\infty (a',s')$ of the
type of node is then given by the solution to the functional equation
\[
f_\infty(a', s')=\int \int f_\infty (a,s)f(a',s'|s,a)dsda.
\]

And, since $f(a',s'|s,a)=f(a'|a)f(s'|s)$ in both versions of the model, it follows that
the age and social index distributions are independent in the
stationary distribution as well: $f_\infty(a', s')=f_\infty
(a')f_\infty (s')$. The corresponding densities are hence the solutions to
\begin{align}
f_\infty (a')&=\int f_\infty (a)f(a'|a)da,\label{f_inf_a}\\
f_\infty (s')&=\int f_\infty (s)f(s'|s)ds. \label{f_inf_s}
\end{align}

The functions $f(a'|a)$ and $f(s'|s)$ were defined in Section 3.1,
$f(a'|a)$ being the same for the two version and $f(s'|s)$ being
different for the U-version and the P-version.
The (unique) solution $f_\infty (a')$ to (\ref{f_inf_a}) is given by
\begin{equation}
f_\infty (a')=\lambda \left(1+\frac{\lambda}{\beta+\mu}\right)e^{-\lambda a'}\left(1-e^{-(\beta+\mu)a'}\right).\label{stat_dist_a}
\end{equation}
As for the solution to (\ref{f_inf_s}), its solution
depends on which version of the model we consider, i.e.\ whether we
use $f^{(P)}(s'|s)$ or $f^{(U)}(s'|s)$.
The (unique) solutions $f^{(P)}_\infty (s')$ and $f^{(U)}_\infty
(s')$  are given by
\begin{equation}
f^{(P)}_\infty (s')=s'f_S(s')/\mu_S, \label{stat_dist_s}
\end{equation}
and
\begin{equation}\label{flU}
f^{(U)}_\infty (s')=\frac{s'+\mu_S}{2\mu_S}f_S(s').
\end{equation}
The P-version hence has the size-biased version of the social index
distribution whereas the U-version has a mixture of the original and
the size-biased social index distribution.

The stationary distribution defined in (\ref{stat_dist_a}) and (\ref{stat_dist_s}) for the P-version, and (\ref{stat_dist_a}) and (\ref{flU}) for the U-version, is hence the distribution of the ''first'' node of a randomly selected edge. Given the type $(a',s')$ of this individual, a randomly selected neighbour of this individual has the previously defined type-distribution $f(a,s|a',s')=f(a|a')f(s|s')$ (where $f(s|s')$, but not $f(a|a')$ depends on which model version we consider).

\subsubsection{The type-distribution obtained using the size-biased method}

We now derive the type distribution  of the nodes of a randomly selected edge using the second method described above, which will be seen to give the same solution. We do this for the P-version of the model.

We know from before (see Equation \ref{mod-uncond-distr}) that the degree $D$ of a randomly selected node in the snap-shot network is mixed Poisson with (random) parameter
\begin{equation}
2\alpha S(1-e^{-(\beta + \mu)A})/(\beta +\mu),\label{rand-par}
\end{equation}
where $S$ has the original social index distribution $f_S$ and $A$ is independent ${\rm Exp}(\lambda)$.
We now instead pick our ''first'' node as a node of a randomly
selected \emph{edge}. As a consequence, this node will have the
size-biased version $\tilde D$ of this random variable (because a node
of degree $k$ has probability proportional to $k$ of being selected). The size-biased distribution of mixed Poisson is also mixed Poisson, but with the size-biased version of the random parameter. And, the size-biased distribution of the random parameter (\ref{rand-par}) equals
\[
\frac{2\alpha}{\beta + \mu}\tilde S\left(\widetilde{1-e^{-(\beta + \mu)A}}\right).
\]
By this tilde-notation we mean that the type $(A^*,S^*)$ of the
''first'' node will have social index distribution $\tilde S$ (the
size biased distribution of the original density, having density
$sf_S(s)/\mu_S$) which hence agrees with (\ref{stat_dist_s}). The age distribution  $A^*$ will be independent and have the distribution defined by $1-e^{-(\beta + \mu)A^*}$ being distributed like $\left(\widetilde{1-e^{-(\beta + \mu)A}}\right)$, where $A\sim {\rm Exp}(\lambda)$. But, since $A\sim {\rm Exp}(\lambda)$ it follows that
\[
1-e^{-(\beta + \mu)A}\sim {\rm Beta}\left( 1,\ \frac{\lambda}{\beta + \mu}\right).
\]
>From this it follows that the size biased distribution
\[
\left(\widetilde{1-e^{-(\beta + \mu)A}}\right) \sim {\rm Beta}\left( 2,\ \frac{\lambda}{\beta + \mu}\right).
\]
So $1-e^{-(\beta + \mu)A^*}$ has this distribution, implying that $e^{-(\beta + \mu)A^*}\sim {\rm Beta}\left( \lambda/(\beta + \mu), 2\right)$. From this it follows after some straightforward calculations that $A^*$ has exactly the density given in Equation (\ref{stat_dist_a}). As in the previous sub-section, given the type $(a',s')$ of the ''first'' node, a randomly selected neighbour of this node has the type-distribution $f(a,s|a',s')=f(a|a')f(s|s')$.

\subsubsection{The degree correlation $\rho$}

In the two previous sub-sections we have derived the type-distribution
of the two nodes of a randomly selected edge: The type $(A^*, S^*)$ of
the ''first'' node has density $f_{\infty}(a')f_{\infty}(s')$ where
$f_{\infty}(a')$ was defined in (\ref{stat_dist_a}) for both versions
of the model, and $f_{\infty}(s')$ was given by (\ref{stat_dist_s})
for the P-version (and not explicit for the U-version). Given the type
$(a',s')$ of the first node, the second node has type-distribution
$f(a,s|a',s')=f(a|a')f(s|s')$, where $f(a|a')$ and $f(s|s')$ were
defined in (\ref{f_U(aa)}) -- (\ref{f_P(ss)}). We now derive the degree correlation between two such randomly selected nodes.

In order to compute the degree correlation of the nodes adjacent to the randomly selected edge we condition on the types. Given the types we know from before that the degree distribution is Poisson: an individual of type $(a,s)$ has a Poisson number of neighbours and the mean $g(a,s)=g_1(a)g_2(s)$ equals
\begin{equation}
g_1(a)= \frac{1-e^{-(\beta + \mu) a}}{\beta + \mu}
\quad\text{and}\quad g_2(s)=
\left\{
\begin{array}{cl}
g^{(U)}_2(s)&=\alpha(s+\mu_S)\qquad\text{for the U-version},\\
g^{(P)}_2(s)&=2\alpha s\qquad\text{for the P-version}.
\end{array}
\right. \label{Po_mean}
\end{equation}

We now compute the degree correlation using these results. Let $D_1$ denote the degree of the ''first'' node and $D_2$ to
its selected neighbour. In order to compute the degree correlation
$\rho (D_1, D_2)$ we need to compute the following expectations:
$\Ex[D_1]=\Ex[D_2]$, $\Ex[D_1^2]=\Ex[D_2^2]$ and $\Ex[D_1D_2]$ (the first as
well as the second moment of $D_1$ and $D_2$ are the same since we have reached stationarity). We do this excluding the edge connecting the two nodes (subtracting the value 1 everywhere does not change the degree correlation). The conditional degree distributions remain Poisson with the same means because if $X\sim {\rm Po}(\beta)$, then the conditional distribution of $X-1$, conditional on that $X\ge 1$, is also ${\rm Po}(\beta)$. We hence get
\begin{align*}
\Ex [D_1]&=\Ex [D_2]=\int g_1(a)f_\infty(a)da\int g_2(s)f_\infty (s)ds,\\
\Ex [D_1^2]&=\Ex [D_2^2]= \int\int (g_1(a)g_2(s)+g_1^2(a)g_2^2(s) )f_\infty(a)f_\infty (s)da ds,\\
\Ex [D_1D_2]&= \int\int g_1(a)g_1(a')f_\infty (a)f(a'|a)dada'   \int\int g_2(s)g_2(s')f_\infty (s)f (s'|s)dsds'.
\end{align*}
Given these expressions, the degree correlation is given by
\begin{equation}
\rho(D_1, D_2)=\frac{\Ex (D_1D_2)-\Ex (D_1)\Ex (D_2)}{\Ex (D_1^2)-(\Ex
  (D_1))^2}= \frac{{\rm C}(D_1,D_2)}{\Var[D_1]}.\label{deg-corr}
\end{equation}
The perhaps most important questions is to learn if $\rho$, or
equivalently the covariance $C$, is positive or negative. For this it is sufficient to compute the numerator of (\ref{deg-corr}).
Using $f_{\infty}(a,s)$ (defined by (\ref{stat_dist_a}),
(\ref{stat_dist_s}), and (\ref{flU}))  and (\ref{Po_mean}) 
in the expressions above for the P-version and U-version, separately,  standard but tedious
calculations reveal that
\begin{align*}
\Ex [D^{(P)}_1]&=\Ex [D^{(P)}_2]=\frac{2}{\lambda +2(\beta+\mu)}\frac{2\alpha \Ex[S^2]}{\Ex[S]},
\\
\Ex [D^{(P)}_1D^{(P)}_2]&=
\frac{5\lambda +6(\beta+\mu)}{(\lambda+\beta+\mu)(\lambda+2(\beta+\mu))(\lambda+3(\beta+\mu))}\left(\frac{2\alpha \Ex[S^2]}{\Ex[S]}\right)^2,
\end{align*}
 as well as
\begin{align*}
\Ex [D^{(U)}_1]&=\Ex
[D^{(U)}_2]=\frac{2}{\lambda+2(\beta+\mu)}\frac{\alpha \left(
  \Ex[S^2]+3(\Ex [S])^2\right)}{2 \Ex[S]},
\\
\Ex [D^{(U)}_1D^{(U)}_2]&=
\frac{5\lambda+6(\beta+\mu)}{(\lambda+\beta+\mu)(\lambda+2(\beta+\mu))(\lambda+3(\beta+\mu))}
2\alpha^2 \left(\Ex[S^2]+(\Ex[S])^2\right).
\end{align*}

This gives us for the
P-version
\begin{equation}\label{cP}
{\rm C}^{(P)}(D_1,D_2)=
\frac{\lambda^2}{(\lambda+\beta+\mu)(\lambda+2(\beta+\mu))^2(\lambda+3(\beta+\mu))} \left(2\alpha \Ex(S^2)\right)^2.
\end{equation}
Since all parameters are positive we conclude that  for the
P-version of the model the covariance, and
hence also the degree correlation $\rho$, is always positive.
This is true irrespective of the model parameters and choices of social index distribution.

For the
U-version the picture is different. Introduce $\gamma=\beta+\mu$ and compute
\begin{align}
{\rm C}^{(U)}(D_1,D_2)&=
\alpha ^2 \frac{6\gamma+5\lambda}{(\lambda+\gamma)(\lambda+2\gamma)(\lambda+3\gamma)}
2 \left(\Ex[S^2]+(\Ex[S])^2\right)
\label{cU}
\\
&\qquad
-\alpha ^2 \frac{1}{(\lambda +
  2\gamma)^2}\left(\frac{\Ex[S^2]+3(\Ex[S])^2}{\Ex[S]}\right)^2\nonumber
\\
&=\frac{\alpha ^2}{(\lambda +
  2\gamma)^2}(\Ex S)^2
\left(\frac{(6\gamma+5\lambda)(\lambda+2\gamma)}{(\lambda+\gamma)(\lambda+3\gamma)}2\left(\frac{\Ex [S^2]}{(\Ex [S])^2}+1 \right)- \left(\frac{\Ex [S^2]}{(\Ex [S])^2}+3 \right)^2\right).
\end{align}
Let us also denote
\[a= \frac{\lambda^2}{(\lambda + \gamma)(\lambda +
  3\gamma)}.\]
Then  we can rewrite (\ref{cU}) as follows
\begin{align}
{\rm C}^{(U)}(D_1,D_2) &=\frac{\alpha ^2}{(\lambda +
  2\gamma)^2}(\Ex S)^2
\left(2(a+4)\left(\frac{\Ex [S^2]}{(\Ex [S])^2}+1 \right)-
  \left(\frac{\Ex [S^2]}{(\Ex [S])^2}+3 \right)^2\right) \label{cU1}
\\
&=\frac{\alpha ^2}{(\lambda +
  2\gamma)^2}(\Ex S)^2
\left(- \left(\frac{\Ex [S^2]}{(\Ex [S])^2}\right)^2 +2\frac{\Ex
    [S^2]}{(\Ex [S])^2}(a+1)-(1-2a)\right).\nonumber
\end{align}
This together with the facts that $\frac{\Ex [S^2]}{(\Ex [S])^2}\geq 1$ and $0<a<1$, yields  that
if
\[\frac{\Ex [S^2]}{(\Ex [S])^2} < 1+a+\sqrt{a^2 +4a}\]
then ${\rm C}^{(U)}(D_1,D_2)>0$, and hence also that $\rho^{(U)}>0$.
If
\[\frac{\Ex [S^2]}{(\Ex [S])^2} = 1+a+\sqrt{a^2 +4a}\]
then
$\rho^{(U)}=0$, and otherwise, $\rho^{(U)}<0$.

We conclude, that for a fixed value of $\Ex [S]$ and the other model parameters, increasing the
variance of $S$ from $0$ (constant $S$) to very large variance, allows us to pass through the
assortative regime, to neutral (no assortativity) and then to
disassortative. One can also observe that increasing the variance
(i.e.\ $\Ex [S^2]$) also increase $R^U$ and $R^P$ respectively (see
(\ref{R^U})) which thus might cross the critical value of 1 leading to
percolation. This second effect, that disassortativity facilitates
percolation, was also reported  by Newman (2002) for a different  model.

\section{Discussion}\label{disc}

In the present paper we studied properties of a stochastic dynamic
network in which nodes are born and die randomly in time, and while
alive pair of nodes are connected and disconnected by edges randomly
in time. The rate of connecting to other nodes depended on the social
index, given at birth, of the node. For this model we derived limiting
properties of a snapshot of the network after a long time, assuming
the node population grew large. 

The three main results were: the type distribution of neighbour nodes
(Section 3.1), a criterion for when phase transition occurs, above
which the network has a giant component (Section 3.2) and an
expression for the degree correlation of connected edges (Section
3.3). One of the most interesting features of the model is that
U-version of the model (in which neighbouring nodes are selected
uniformly among all living nodes) may have either positive or negative
degree correlation depending on the numerical values of the model
parameters (birth and death rates of nodes and edge) and the social
index distribution. 

Together with previous results of the model (Britton and Lindholm,
(2010)) the most important limiting local properties of the network
are hence known by now, as well as whether or not the network has a
giant component. Other global properties such as the diameter of the
network remains to be analysed. Another interesting class of problems
to study for the present model would be to study dynamic properties of
the network. One such problem would be to study limiting properties of
some process taking place ''on'' the network. For example, can an
epidemic, forest fire or similar persist on the network forever or
will it die out? What happens if a lightning process kills randomly
selected nodes and all their neighbours at a constant rate. 

The model can of course also be generalized in several ways to
make it more applicable for certain situations. Nodes could be of
different types with different attachment rates depending on the types
of the nodes in question, for example male or female if mimicking a sexual
network. Similarly, edges could be of different types, 
perhaps reflecting the ''degree'' of acquaintance which for example
could affect the transmission probability of an epidemic taking place
on the network. Still, the current model does capture some important
properties of a real world network in the nodes as well as edges are
created and cease to exist. Also, depending on the model parameters,
the model is quite flexible in producing different degree
distributions and degree correlations (the network has no clustering
asymptotically). 

\section*{Acknowledgments}

T.B.\ is grateful to Riksbankens jubileumsfond (The Bank of
Sweden Tercentenary Foundation). Part of this work was done while the
authors were visiting Institut Mittag-Leffler to which we are
grateful. We also thank Pieter Trapman for fruitful discussions.

\end{document}